\documentstyle[11pt]{article}
\title{On arrangements of the roots of a hyperbolic polynomial and of one of 
its derivatives}
\author{Vladimir Petrov Kostov\\ \hspace{8cm}{\sl To Prof. Rumyan Lazov}} 
\date{}
\newtheorem{tm}{Theorem}
\newtheorem{defi}[tm]{Definition}
\newtheorem{rem}[tm]{Remark}
\newtheorem{rems}[tm]{Remarks}
\newtheorem{lm}[tm]{Lemma}
\newtheorem{ex}[tm]{Example}
\newtheorem{cor}[tm]{Corollary}
\newtheorem{prop}[tm]{Proposition}
\newtheorem{nota}[tm]{Notation}
\begin{document}
\maketitle 

\begin{abstract}
We consider real monic 
{\em hyperbolic} polynomials in 
one real variable, i.e. polynomials having only real roots.
Call {\em hyperbolicity domain} $\Pi$ of the family of polynomials 
$P(x,a)=x^n+a_1x^{n-1}+\ldots +a_n$, 
$a_i,x\in {\bf R}$, the set 
$\{ a\in {\bf R}^n| P~$is hyperbolic~$\}$. The paper studies a stratification 
of $\Pi$ defined by the arrangement of the roots of $P$ and $P^{(k)}$, 
where $2\leq k\leq n-1$. We prove that the strata are smooth contractible 
real algebraic varieties.\\


\end{abstract}

\section{Introduction}

In the present paper we consider real monic 
{\em hyperbolic} (resp. {\em strictly hyperbolic)} polynomials in 
one real variable, i.e. polynomials having only real (resp. only 
real distinct) roots. If a polynomial is (strictly) hyperbolic, then such 
are all its non-trivial derivatives.

Consider the family of polynomials $P(x,a)=x^n+a_1x^{n-1}+\ldots +a_n$, 
$a_i,x\in {\bf R}$. Call {\em hyperbolicity domain} $\Pi$ the set 
$\{ a\in {\bf R}^n| P~$is hyperbolic~$\}$. The paper studies a stratification 
of $\Pi$ defined by the {\em configuration} (we write sometimes 
{\em arrangement}) of the roots of $P$ and $P^{(k)}$, 
where $2\leq k\leq n-1$. The study of this stratification 
began in \cite{KoSh}, see also \cite{Ko1} and \cite{Ko2} for the particular 
cases $n=4$ and $n=5$. Properties of $\Pi$ were proved in 
\cite{Ko3} and \cite{Ko4}, the latter two papers use results 
of V.I. Arnold (see \cite{Ar}), A.B. Givental (see \cite{Gi}) and 
I. Meguerditchian (see \cite{Me1} and \cite{Me2}). 

\begin{nota}
Denote by $x_1\leq \ldots \leq x_n$ the roots of $P$ and by 
$\xi _1\leq \ldots \leq \xi _{n-k}$ the ones of $P^{(k)}$. We write sometimes 
$x_i^{(k)}$ instead of $\xi _i$ if the index $k$ varies. Denote by 
$y_1<\ldots <y_q$ the {\em distinct} roots of $P$ and by $m_1$, $\ldots$, 
$m_q$ their multiplicities (hence, $m_1+\ldots +m_q=n$).
\end{nota} 

The classical 
Rolle theorem implies that one has the following chain of inequalities: 

\begin{equation}\label{Rolle}
x_i\leq \xi _i\leq x_{i+k}~~,~~i=1,\ldots ,n-k
\end{equation} 

\begin{defi}\label{CV}
A {\em configuration vector (CV)} of {\em length} $n$
is a vector whose components are either
positive integers (sometimes indexed by the letter $a$, their sum
being $n$) or the letter $a$. The integers equal the multiplicities of the
roots of $P$, the letters $a$ indicate the positions of the
roots of $P^{(k)}$; $m_a$ means that a root of $P$ of multiplicity 
$m<k$ coincides with a simple root of $P^{(k)}$.
A CV is called {\em a priori admissible} if for the configuration of the
roots of $P$ and $P^{(k)}$ defined by it there hold 
inequalities (\ref{Rolle}).
\end{defi}

\begin{rem}\label{simpleroot}
If a root of $P$ of multiplicity $<k$ is also a root of $P^{(k)}$, 
then it is a simple root of $P^{(k)}$, see Lemma 4.2 from \cite{KoSh}. By 
definition ``a root of multiplicity $0$'' means ``a non-root''.
\end{rem}
 
\begin{ex}
For $n=8$, $k=3$ the CV $(1,a,1,2_a,a,a,4)$ (which is a priori admissible) 
means that the roots $x_j$ and
$\xi _i$ are situated as follows:
$x_1<\xi _1<x_2<x_3=x_4=\xi _2<\xi _3<\xi _4<x_5=\ldots =x_8=\xi _5$. The 
multiplicity $4$ is not indexed with $a$ because it is $>k$, i.e. it 
automatically implies $x_5=\ldots =x_8=\xi _5$.
\end{ex}

\begin{defi}\label{classAB}
Given a hyperbolic polynomial $P$ call {\em roots of class B} (resp. 
{\em roots of class A}) its roots of multiplicity $<k$ which 
coincide with roots of $P^{(k)}$ (resp. all its other roots). In a CV the 
roots of class B correspond to multiplicities indexed by $a$.
\end{defi}

\begin{defi} 
For a given CV $\vec{v}$ call {\em stratum} of $\Pi$ (defined by $\vec{v}$) 
its subset of polynomials $P$ 
with configuration of the roots of $P$ and $P^{(k)}$ defined by $\vec{v}$.
\end{defi}

The aim of the present paper is to prove the following 

\begin{tm}\label{stratificationtm}
All strata of this stratification are smooth 
contractible real algebraic varieties. 
\end{tm}

The theorem is proved in Section~\ref{prooftm}.

\begin{rem}
It is shown in \cite{KoSh}, Theorem 4.4, that every a priori admissible CV 
defines a non-empty connected stratum. The essentially new result of the 
present paper is the proof not only of connectedness but of contractibility. 
In \cite{Ko5} the notion of a priori admissible 
CV is generalized in the case of not necessarily hyperbolic polynomials and it 
is shown there that all such CVs are realizable by the arrangements of the 
real roots of polynomials $P$ and of their derivatives $P^{(k)}$ (the position 
and multiplicity of the complex roots is not taken into account there).
\end{rem}

\begin{nota}
We denote by $D(i,j)$ the {\em discriminant set} 
$\{ a\in {\bf R}^n|{\rm Res}(P^{(i)},P^{(j)})=0\}$ (recall that for 
$a\in \Pi$ one has 
{\rm Res}$(P^{(i)},P^{(j)})=0$ if and only if $P^{(i)}$ and $P^{(j)}$ have 
a common root).
\end{nota}

Denote by $G$ a point from $\Pi$. Consider the discriminant set $D(0,k)$, 
$k\geq 2$, at $G$ for $G$ lying strictly 
inside $\Pi$ at which there hold exactly $s$ equalities of the form 
$x^{(k)}_j=x_i$, with $s$ different indices $j$ and $s$ different 
indices $i$.

\begin{prop}\label{stratificationprop}
In a neighbourhood of the point $G$ the set $D(0,k)$ is locally the 
union of $s$ smooth hypersurfaces intersecting transversally at $G$. 
\end{prop}

All propositions are proved in Section~\ref{proofprop}. 
The proposition can be generalized in the following way. 
Suppose that at a point $G$ lying strictly inside $\Pi$ there hold exactly 
$s$ equalities $x_j^{(k_i)}=x_i$, with $s$ different indices $i$ and $s$ 
different couples $(k_i,j)$. 

\begin{prop}\label{stratificationprop1}
In a neighbourhood of the point $G$ these $s$ equalitites define $s$ 
smooth hypersurfaces intersecting transversally at $G$. 
\end{prop}

\begin{rem}\label{MultV}
It is shown in \cite{Ko3} that for each $q$-tuple of positive integers 
$m_j$ with sum $n$ the subset $T$ of $\Pi$ (we call it a {\em stratum} 
of $\Pi$ defined by the {\em multiplicity vector} $(m_1,\ldots ,m_q)$, 
not by a CV) consisting of 
polynomials with distinct roots $y_i$, 
of multiplicities $m_i$, is a smooth variety of dimension 
$q$ in ${\bf R}^n$. 
\end{rem}

Denote by $T$ a stratum of $\Pi$ defined by a multiplicity vector. 
Fix a point $G\in T$. 
Suppose that at $G$ there are $s$ among the roots $y_j$ which are of class B.
Suppose that one has $m_i<k$ for all $i$. 
The condition $m_i<k$ implies that all points from $D(0,k)\cap T$ close to $G$ 
result from roots of $P^{(k)}$ coinciding with roots of $P$ of class B.

\begin{prop}\label{stratificationpropbis}
In a neighbourhood of the point $G$ the set $D(0,k)\cap T$ is locally the 
union of $s$ smooth codimension $1$ subvarieties of $T$ intersecting 
transversally at $G$.
\end{prop}

\begin{rems}
1) A stratum of $\Pi$ of codimension $\kappa \leq k$ defined by $\kappa$ 
equalities of the 
form $x_i=\xi _j$ (i.e. $P$ has no multiple root) has a tangent space 
transversal to the space $Oa_{n-\kappa +1}\ldots a_n$. Indeed, the roots 
$\xi _j$ depend smoothly on $a_1,\ldots ,a_{n-k}$, and the conditions 
$P(\xi _j,a)=0$ allow one to express $a_{n-\kappa +1},\ldots ,a_n$ as smooth 
functions of $a_1,\ldots ,a_{n-\kappa}$ (use Vandermonde's determinant with 
distinct arguments $\xi _1$, $\ldots$, $\xi _{\kappa}$). 
It would be nice to prove or disprove the 
statements: 

A) this property holds without the assumption $\kappa \leq k$ and 
that $P$ has no multiple root; 

B) the limit of the tangent space to the stratum when a 
stratum of lower dimension from its closure is approached exists and is 
transversal to the space $Oa_{n-\kappa +1}\ldots a_n$. 

For $n=4$ and $n=5$ this seems to be  
true, see \cite{Ko1} and \cite{Ko2}. 
The statements would be a generalization of such a transversality 
property of the strata of $\Pi$ defined by multiplicity vectors, not by 
CVs (proved in \cite{Ko3}, Theorem~1.8; see Remark~\ref{MultV}). 
Outside $\Pi$ the first statement is not true -- for $n=4$, $a_1=0$, 
the discriminant set $D(0,2)$ has a Whitney umbrella singularity at the origin 
and there are points where its tangent space is parallel to $Oa_4$; this 
can be deduced from \cite{Ko1} (see Section~$3$ and Lemma~29 in it). 

2) In \cite{KoSh}, \cite{Ko1} and \cite{Ko2} a stratification of $\Pi$ defined 
by the arrangement of all 
roots of $P$, $P'$, $\ldots$, $P^{(n-1)}$ is considered (the initial idea to 
consider this stratification belongs to B.Z. Shapiro). The results 
of the present paper cannot be transferred directly to that case for two 
reasons: 

a) for $n\geq 4$ not all arrangements consistent with (\ref{Rolle}) are 
realized by 
hyperbolic polynomials and it is not clear how to determine for any 
$n\in {\bf N}^*$ the realizable 
ones (e.g. for $n=4$ only $10$ out of $12$ such 
arrangements are realized, see \cite{KoSh} or \cite{Ko1}; 
for $n=5$ only $116$ out of $286$, see \cite{Ko1}); 
the reason for this is clear -- a 
monic polynomial has only $n$ coefficients that can be varied 
whereas there are $n(n+1)/2$ roots of $P$, $P'$, $\ldots$, $P^{(n-1)}$;

b) for $n\geq 4$ there are {\em overdetermined strata}, 
i.e. strata on which the number 
of equalities between any two of the roots of 
$P$, $P'$, $\ldots$, $P^{(n-1)}$ is greater than the codimension of the 
stratum.  
\end{rems} 

In Section~\ref{techlm} we 
prove two technical lemmas (and their corollaries) 
used in the proof of the theorem and the propositions. 
Section~\ref{dimensionCV} 
is devoted to the dimension of a stratum and its relationship with the CV 
defining it. The above  
propositions are just the first steps in the study of the 
set $D(0,1)\cup D(0,k)$ (and, more generally, of the set 
$D(0,1)\cup \ldots \cup D(0,n-1)$) at a point of $\Pi$. 

\section{Configuration vectors and dimensions of strata
\protect\label{dimensionCV}}

In this section we recall briefly results some of which are 
from \cite{KoSh}:\\  

1) Call {\em excess of multiplicity} of a CV  
the sum $m=\sum (m_j-1)$ taken over all multiplicities $m_j$ of distinct 
roots of $P$. 
A {\em stratum} of codimension
$i$ is defined by a CV which has exactly $i-m$ letters $a$ as indices, i.e. 
the polynomial $P$ has exactly $i-m$ distinct roots of class B.\\ 

2) A stratum of codimension $i$ is locally a smooth real algebraic variety of
dimension $n-i$ in ${\bf R}^n$.\\ 

3) In what follows we say a stratum of codimension $i$ to be of dimension 
$n-i-2$. We decrease its dimension in ${\bf R}^n$ by $2$ 
to factor out the possible shifting of the 
variable $x$ by a constant and the one-parameter group of transformations 
$x\mapsto \exp (t)x$, $a_j\mapsto \exp (jt)a_j$, $t\in {\bf R}$; 
both of them leave CVs unchanged. This allows one to consider the family 
$P$ only for $a_1=0$, $a_2=-1$ (if $a_1=0$, then 
there are no hyperbolic polynomials 
for $a_2>0$ and for $a_2=0$ the only one is $x^n$).\\ 

4) In accordance with the convention from 3), it can be deduced from 1) that 
the CVs defining strata of 
dimension $\delta$ are exactly the ones in which the polynomial $P$ has 
$\delta +2$ distinct roots of 
class A, i.e. these are CVs having 
$\delta +2$ components which are multiplicities of roots 
of $P$ not indexed by the letter $a$.\\

5) A point of a stratum of codimension $i>1$ defined by a CV $\vec{v}$
belongs to the closure of any stratum of codimension $i-1$
whose CV $\vec{w}$ is
obtained from $\vec{v}$ by means of one of the following three operations:

i) if $\vec{v}=(A,l_a,B)$, $l\leq k-1$, $A$ and $B$ are non-void,
then $\vec{w}=(A,l,a,B)$ or $\vec{w}=(A,a,l,B)$;

ii) if $\vec{v}=(A,r_a,B)$, $r\leq k-1$, $A$ and $B$ are non-void,
then $\vec{w}=(A,r',r''_a,B)$ or $\vec{w}=(A,r'_a,r'',B)$, $r'>0$, $r''>0$,
$r'+r''=r$.

iii) if $\vec{v}=(A,r,B)$,
then $\vec{w}=(A,C,B)$ where $C$ is a CV defining a stratum of dimension 
$0$ in ${\bf R}^r$, see 4).\\ 

6) It follows from the definition of the codimension of a
stratum that the
three possibilities {\em i)}, {\em ii)} and {\em iii)} from 5) 
are the only ones
to increase by 1 the dimension of a stratum $S$ when passing to a stratum
containing $S$ in its closure. Indeed, one has to increase by $1$ the number 
of roots of class A, see 4). 
If to this end one has to change the number or the 
multiplicities of the roots 
of class B, then there are no possibilities other than {\em i)} and {\em ii)}. 
If not, then exactly one root $x_i$ of class A must bifurcate, the roots 
stemming from it and the roots of $P^{(k)}$ close to $x_i$ must define 
an a priori admissible CV (they must satisfy conditions (\ref{Rolle})),  
and among these roots  
there must be exactly two of class A. Hence, the bifurcating roots must 
define a CV of dimension $0$ in ${\bf R}^r$, see 4).

\section{Two technical lemmas and their corollaries\protect\label{techlm}}

For a monic strictly hyperbolic polynomial $P$ of degree $n$ consider the 
roots $x^{(k)}_j$ of $P^{(k)}$ as functions of the roots $x_i$ of $P$. Hence, 
these functions are smooth because the roots $x_j^{(k)}$ are simple, see 
Remark~\ref{simpleroot}.

\begin{lm}\label{positive}
For $i=1,\ldots ,n$, $k=1,\ldots ,n-1$, $j=1,\ldots ,n-k$ one has 
$\frac{\partial (x^{(k)}_j)}{\partial (x_i)}>0$.
\end{lm}

{\em Proof:}

$1^0$. Set $x_i=c$, $P=(x-c)Q(x)$, deg$Q=n-1$, $\xi =x^{(k)}_j$. Prove that 
for $k=1$ one has 
$\frac{\partial \xi }{\partial c}>0$. One has 
$(\xi -c)Q'(\xi )+Q(\xi )=0$. Hence, 

\[ \left( \frac{\partial \xi }{\partial c}-1\right) Q'(\xi )+
(\xi -c)Q''(\xi )\frac{\partial \xi }{\partial c}+
Q'(\xi )\frac{\partial \xi }{\partial c}=0~,~{\rm i.e.}\]

\[ \frac{\partial \xi }{\partial c}=\frac{Q'(\xi )}
{(\xi -c)Q''(\xi )+2Q'(\xi )}=
\frac{Q'(\xi )}{P''(\xi )}~.\]

As $Q(\xi )=\frac{P(\xi )}{\xi -c}$, one has 

\begin{equation}\label{deriv}
\frac{\partial \xi }{\partial c}=\frac{(\xi -c)P'(\xi )-P(\xi )}
{(\xi -c)^2P''(\xi )}=
-\frac{P(\xi )}{(\xi -c)^2P''(\xi )}
\end{equation}

For a strictly hyperbolic monic polynomial the signs of $P(\xi )$ and 
$P''(\xi )$ 
are opposite and $\xi \neq c$. This proves the lemma for $k=1$.\\ 

$2^0$. For $k>1$ use induction on $k$. Considering the roots of $P^{(k+1)}$ as 
functions of the ones of $P^{(k)}$ one can write

\begin{equation}\label{complex} 
\frac{\partial (x^{(k+1)}_j)}{\partial c}=\sum _{\nu =1}^{n-k}
\frac{\partial (x^{(k+1)}_j)}{\partial (x^{(k)}_{\nu})}\, 
\frac{\partial (x^{(k)}_{\nu})}{\partial c}
\end{equation} 
and observe that all factors in the right hand-side are $>0$. The lemma is 
proved.~~~~$\Box$

\begin{rem}
The roots $x_j^{(k)}$ are $C^1$-smooth functions of the roots $x_i$ (one can 
forget for a moment that $x_1\leq \ldots \leq x_n$ and assume that 
$(x_1,\ldots ,x_n)\in {\bf R}^n$ and the claim is true for not necessarily 
strictly hyperbolic polynomials; however, in order to define correctly 
$x_j^{(k)}$ one has 
to impose the condition $x_1^{(k)}\leq \ldots \leq x_{n-k}^{(k)}$). Indeed, 
it suffices to prove this for $k=1$ (because  
in the same way one proves that the roots 
of $P^{(\nu +1)}$ are $C^1$-smooth functions of the roots of $P^{(\nu )}$ for 
$\nu =1,\ldots ,n-2$ etc.). For $k=1$ the claim can be deduced from equality 
(\ref{deriv}) -- the fraction in the right hand-side 
has a finite limit for $\xi \rightarrow c$ 
(this limit depends on the order of $c$ as a zero of $P$) and 
for $\xi$ close to $c$ it is a function continuous in $c$. We leave the 
details for the reader. 
\end{rem}  

\begin{cor}\label{positivecor}
For a (not necessarily strictly) hyperbolic polynomial one has 
$\frac{\partial (x^{(k)}_j)}{\partial (x_i)}\geq 0$ for $i,j,k$ like in the 
lemma.
\end{cor}

The corollary is automatic.

\begin{cor}\label{estim}
For a monic strictly hyperbolic polynomial one has 
$0<\frac{\partial (x^{(k)}_j)}{\partial (x_i)}<\frac{n-k}{n}$ 
(for $i,j,k$ like in the lemma) and 
$\sum _{j=1}^{n-k}\frac{\partial (x^{(k)}_j)}{\partial (x_i)}=\frac{n-k}{n}$. 
\end{cor}

{\em Proof:}

By Vieta's formulas one has $x_1+\ldots +x_n=-a_1$, 
$x^{(k)}_1+\ldots +x^{(k)}_{n-k}=-\frac{n-k}{n}a_1$. As 
$\frac{\partial (x^{(k)}_j)}{\partial (x_i)}>0$ for all $j$, one has 

\[ \frac{\partial (x^{(k)}_j)}{\partial (x_i)}<
\frac{\partial (x^{(k)}_1+\ldots +x^{(k)}_{n-k})}{\partial (x_i)}=
\frac{(n-k)}{n}
\, \frac{\partial (x_1+\ldots +x_n)}{\partial (x_i)}=\frac{n-k}{n}\]
which proves the corollary.~~~~$\Box$

\begin{rem}
In the above corollary one sums up w.r.t. the index $j$. When summing up 
w.r.t. $i$ one obtains the equality 
\begin{equation}\label{sumi}
\sum _{i=1}^n\frac{\partial (x^{(k)}_j)}{\partial (x_i)}=1
\end{equation} 
Indeed, if the 
roots $x_i$ are functions of one real parameter (say, $\tau$), 
then one has the equality 
$\sum _{i=1}^n\frac{\partial (x^{(k)}_j)}{\partial (x_i)}\dot{x}_i=
\dot{x}^{(k)}_j$ where $\dot{x}_i$ stands for $\frac{dx_i}{d\tau}$. 
When one has $\dot{x}_i=1$ for all $i$, i.e. the variable 
$x$ is shifted with constant speed $1$, then one has $\dot{x}^{(k)}_j=1$ for 
all $k,j$ and one gets (\ref{sumi}). One needs not suppose the roots $x_i$ 
distinct.
\end{rem} 

In the case of a not strictly hyperbolic polynomial $P$ consider the roots 
$x_j^{(k)}$ as functions of the distinct roots $y_i$ of $P$ (their 
multiplicities remain fixed). 

\begin{lm}\label{positive1}
For $i=1,\ldots ,q$, $k=1,\ldots ,n-1$ one has   
$\frac{\partial (x_j^{(k)})}{\partial (y_i)}\geq 0$ with equality exactly if 
$x_j^{(k)}$ is a root of $P$ of multiplicity $\geq k$ (hence, of multiplicity 
$\geq k+1$) and $x_j^{(k)}\neq y_i$.
\end{lm}

{\em Proof:}

$1^0$. The proof follows the same ideas as the proof of Lemma~\ref{positive}. 
Set $\xi =x_j^{(k)}$, $P=(x-c)^sQ(x)$ where $c=y_i$, $s=m_i$ for some 
$i$, $1\leq i\leq q$. 

Let first $k=1$. One has 
$(\xi -c)^sQ'(\xi )+s(\xi -c)^{s-1}Q(\xi )=0$. Hence, 

\[ s\left( \frac{\partial \xi }{\partial c}-1\right) (\xi -c)^{s-1}Q'(\xi )+
(\xi -c)^sQ''(\xi )\frac{\partial \xi }{\partial c}+\]
\[ +s(\xi -c)^{s-1}Q'(\xi )\frac{\partial \xi }{\partial c}+
s(s-1)(\xi -c)^{s-2}\left( \frac{\partial \xi }{\partial c}-1\right) 
Q(\xi )=0~,~
{\rm i.e.}\]

\[ \frac{\partial \xi }{\partial c}=\frac{s(s-1)Q(\xi )+s(\xi -c)Q'(\xi )}
{s(s-1)Q(\xi )+2s(\xi -c)Q'(\xi )+(\xi -c)^2Q''(\xi )}=\]
\[ =\frac{s((\xi -c)P'(\xi)-P(\xi ))}{(\xi -c)^2P''(\xi )}~.\]
If $\xi =c$, i.e. $s>1$, then $\frac{\partial \xi }{\partial c}=1$. 
If not, then $\frac{\partial \xi }{\partial c}=-\frac{P(\xi )}
{(\xi -c)^2P''(\xi )}$. 
Either $P(\xi )=P'(\xi )=0$ and in this case 
$\frac{\partial \xi }{\partial c}=0$ whatever the multiplicity of $\xi$ as a 
root of $P$ is, or $P(\xi )\neq 0$, 
$P(\xi )$ and $P''(\xi )$ have opposite signs and 
$\frac{\partial \xi }{\partial c}>0$. This proves the lemma for $k=1$.\\ 

$2^0$. For $k>1$ use induction on $k$. Consider the roots of $P^{(k+1)}$ as 
functions of the ones of $P^{(k)}$. Then there holds (\ref{complex}).  
All factors in the right hand-side are $\geq 0$. 

One has $\frac{\partial (x^{(k+1)}_j)}{\partial c}=0$ exactly if 
in every summand in the right hand-side of (\ref{complex}) 
at least one of the two factors is $0$. This is the case if 
$\xi =x^{(k+1)}_j$ is a root of $P$ of multiplicity $\geq k+1$ and 
$\xi \neq c$. 
Indeed, in this case one has 
$\frac{\partial (x^{(k+1)}_j)}{\partial (x^{(k)}_{\nu})}=0$ if  
$x^{(k+1)}_j\neq x^{(k)}_{\nu}$ and 
$\frac{\partial (x^{(k)}_{\nu})}{\partial c}=0$ if $x^{(k+1)}_j=x^{(k)}_{\nu}$ 
(and, hence, $x^{(k)}_{\nu}\neq c$). 

If $\xi$ is a root of $P$ of multiplicity $\geq k+1$ and 
$\xi =c$, then one has $\frac{\partial (x^{(k+1)}_j)}{\partial c}=1$.

If $\xi$ is a root of $P$ of multiplicity $\leq k$, then it is not a root 
of $P^{(k)}$. Hence, 
$\frac{\partial (x^{(k+1)}_j)}{\partial (x^{(k)}_{\nu})}>0$ 
for all $\nu$. At least one of the factors 
$\frac{\partial (x^{(k)}_{\nu})}{\partial c}$ is $>0$ 
(i.e. for at least one $\nu$). 
Indeed, if $c$ is a root of $P$ of multiplicity $\geq k+1$, then this is true 
for the root $x^{(k)}_{\nu}$ which equals $c$ (by inductive assumption). 
If $c$ is a root of $P$ of 
multiplicity $\leq k$, then there exists a simple root $x^{(k)}_{\nu}$ 
of $P^{(k)}$ (this follows from Rolle's theorem applied $k$ times). 
Hence, $x^{(k)}_{\nu}$ is a root of $P$ of multiplicity $\leq k-1$, and  
for this root one has $\frac{\partial (x^{(k)}_{\nu})}{\partial c}>0$.

The lemma is proved.~~~~$\Box$ 

\begin{cor}\label{estimbis}
For a monic hyperbolic polynomial one has 
$0\leq \frac{\partial (x^{(k)}_j)}{\partial (y_i)}\leq \frac{n-k}{n}$ 
(for $i,j,k$ like in the lemma) and 
$\sum _{j=1}^{n-k}\frac{\partial (x^{(k)}_j)}{\partial (y_i)}=\frac{n-k}{n}$. 
\end{cor}

The corollary is proved by analogy with Corollary~\ref{estim}. 

\section{Proofs of the propositions\protect\label{proofprop}}

{\bf Proof of Proposition~\ref{stratificationprop}:}\\ 

Prove the smoothness. The roots $x^{(k)}_j$ are smooth functions of the 
coefficients $a_1$, $\ldots$, $a_{n-k}$. The condition $P(x^{(k)}_j,a)=0$ 
allows one to express $a_n$ as a smooth function of $a_1$, $\ldots$, 
$a_{n-1}$. Hence, this equation defines locally a smooth hypersurface in 
${\bf R}^n$.  

To prove the transversality assume first that the indices are changed so that 
$i=j=1,\ldots ,s$. It suffices to prove that the ``Jacobian'' matrix 
$\left\{ \frac{\partial (x_j-x^{(k)}_j)}{\partial (x_{\nu })}\right\}$, 
$j,\nu =1,\ldots ,s$, 
is of maximal rank (in the true Jacobian matrix one has 
$\nu =1,\ldots ,n$, not $\nu =1,\ldots s$). Its diagonal entries equal 
$1-\frac{\partial (x^{(k)}_j)}{\partial (x_j)}$ while its non-diagonal ones 
equal 
$-\frac{\partial (x^{(k)}_j)}{\partial (x_\nu )}$. 
Corollary~\ref{estim} implies 
that the matrix is diagonally dominated -- for $\nu$ fixed its diagonal entry 
(which is positive) is greater than the sum of the absolute values of its 
non-diagonal entries (which are all negative). Hence, the matrix is 
non-degenerate.~~~~$\Box$\\ 

{\bf Proof of Proposition~\ref{stratificationprop1}:}\\

The proof of the smoothness is done like in the proof of 
Proposition~\ref{stratificationprop}. To prove the transversality assume again 
that $i=j=1,\ldots ,s$ and consider again 
the ``Jacobian'' matrix 
$\left\{ \frac{\partial (x_j-x^{(k_j)}_j)}{\partial (x_{\nu })}\right\}$, 
$j,\nu =1,\ldots ,s$. 
Like in the previous proof we show that the matrix is diagonally dominated, 
hence, non-degenerate. However, the numbers $k_j$ are not necessarily the 
same and therefore we fix $j$ (hence, $k_j$ as well) and we change $\nu$. 
By equality (\ref{sumi}), one has 
  
\[ \sum _{\nu =1}^s\frac{\partial (x_j-x^{(k_j)}_j)}{\partial (x_{\nu })}=
1-\sum _{\nu =1}^s\frac{\partial (x^{(k_j)}_j)}{\partial (x_{\nu })}\geq 
1-\sum _{\nu =1}^n\frac{\partial (x^{(k_j)}_j)}{\partial (x_{\nu })}=0\]
and the case of equality has to be excluded because the smallest and the 
greatest root of $P$ are not 
among the roots $x_1$, $\ldots$, $x_s$ and all partial derivatives are 
strictly positive, see Lemma~\ref{positive}. The last inequality implies that 
the matrix is diagonally dominated. ~~~~$\Box$\\

{\bf Proof of Proposition~\ref{stratificationpropbis}:}\\ 

The proof is almost a repetition of the one of 
Proposition~\ref{stratificationprop}. The only difference is that the 
Jacobian matrix looks like this: 
$\left\{ \frac{\partial (y_j-m_{\nu }x^{(k)}_j)}{\partial (y_{\nu })}\right\}$ 
(recall that $y_{\nu}$, of multiplicity $m_{\nu}$, are the 
distinct roots of $P$).~~~~$\Box$

\section{Proof of Theorem \protect\ref{stratificationtm}
\protect\label{prooftm}}

$1^0$. Smoothness is proved in \cite{KoSh}, Proposition 4.5; algebraicity is 
evident. So one has to 
prove only contractibility. Assume that $a_1=0$, $a_2=-1$. 

To prove contractibility of the strata represent each 
stratum $T$ of dimension $\delta \geq 1$ 
as a fibration  
whose fibres are one-dimensional varieties with 
the following properties: 

a) the fibres are phase curves of a smooth 
vectorfield without stationary points defined on $T$; hence, each fibre 
can be smoothly parametrized by $\tau \in (0,1)$; 
this is proved in $2^0$ -- $4^0$;  

b) the limits for $\tau \rightarrow 1$ of the points of 
the fibres exist and they belong to a finite union ${\cal U}$ 
of strata of lower 
dimension; we call the limits {\em endpoints}; the proof of this is given in 
$3^0$ -- $5^0$;

c) the union ${\cal U}$ is a contractible set (proved in $7^0$ -- $8^0$);

d) each point of the union ${\cal U}$ is the endpoint of some fibre (proved in 
$6^0$).

Thus the union ${\cal U}$ is a retract of the given stratum and 
contractibility of ${\cal U}$ implies the one of the stratum. Contractibility 
of the strata of dimension $0$ will be proved directly (in $7^0$).\\

$2^0$. A shift $\gamma _1$ and a rescaling $\gamma _2$ 
of the $x$-axis fix the smallest root of $P$ at $0$ 
and the greatest one at $1$. Set 
$\gamma =\gamma _2\circ \gamma _1$.

\begin{nota}
Denote by $\Delta$ the set of monic hyperbolic polynomials 
obtained from the stratum $T$ by applying the transformation $\gamma$ to each 
point of $T$. 
\end{nota}

\begin{rem}
The set $\Delta$ (like $T$) is 
a smooth variety of dimension $\delta$. The transformation $\gamma$  
defines a diffeomorphism $\bar{T}\rightarrow \bar{\Delta}$ while 
$\gamma ^{-1}$ defines  
a diffeomorphism $\bar{\Delta} \rightarrow \bar{T}$; this can be deduced 
from the conditions $a_1=0$, $a_2=-1$. 
\end{rem}

$3^0$. Recall that $y_i$ denote the distinct roots of $P$. 
We construct (see $4^0$ -- $5^0$) the speeds $\dot{y}_i$ on $\Delta$ 
which amounts to constructing a 
vectorfield defined on $\Delta$. Therefore the fibration from $1^0$ can 
be defined by means of the phase curves of a vectorfield defined 
on $T$ (to this end one has to apply $\gamma ^{-1}$). 
We leave the technical details for the reader. 

\begin{rem}
It follows from our 
construction (see in particular part 3) of Lemma~\ref{speed}) that these 
two vectorfields can be continuously extended respectively on $\bar{\Delta}$ 
and $\bar{T}$.
\end{rem}

Along a phase curve of the vectorfield, all roots of $P$ of 
class A except one (in particular, the smallest and the greatest one) 
do not change their position and 
multiplicity; the rest of the roots of $P$ do not 
change their multiplicity. The limits (forwards and 
backwards) of the points of the phase curves 
exist when the boundary of $\Delta$ is approached. At these 
limit points, if a confluence of roots of $P$ occurs, then the 
multiplicities of the coinciding roots are added. The images 
under $\gamma ^{-1}$ of the forward limits are 
the endpoints (see b) from $1^0$).

Denote by $P_{\sigma }$ ($\sigma \in {\bf R}$) a family of monic 
hyperbolic polynomials 
represented by the points of a given phase curve in $\Delta$. 
We prove in $4^0$ that there exists $\sigma _0>0$ 
such that for $\sigma \in [0,\sigma _0)$ one has $P_{\sigma }\in \Delta$ 
(hence, $\gamma ^{-1}(P_{\sigma })\in T$) while $P_{\sigma _0}\not\in \Delta$ 
(hence, $\gamma ^{-1}(P_{\sigma _0})\not\in T$). The polynomial 
$P_{\sigma _0}$ represents the forward limit point of the 
given phase curve. We set 
$\dot{y}_i=dy_i/d\sigma$.\\ 

$4^0$. Change for convenience (in $4^0$ -- $6^0$) the indices of the distinct 
roots $y_i$ of $P$ and of the roots  
$\xi _i$ of $P^{(k)}$. Choose a root of class A different from the 
smallest and the 
greatest one. Denote it by $y_1$. Denote by $y_2$, $\ldots$, 
$y_d$ the roots of class B and by $\xi _2$, $\ldots$, $\xi _d$ the roots of 
$P^{(k)}$ which are equal to them.

Set $\dot{y}_1=1$. We  
look for speeds $\dot{y}_i$ for which one has $\dot{y}_i=\dot{\xi }_i$, 
$i=2,\ldots ,d$. Hence, one would have 
$y_i=\xi _i$, $i=2,\ldots ,d$, and the multiplicities of the roots of 
$P$ do not change for $\sigma >0$ close to $0$. 
This means that for all such values of $\sigma$ for which the 
order of the union of roots of $P$ and $P^{(k)}$ is preserved, the point 
$\gamma ^{-1}(P_{\sigma})$ belongs to $T$. The value  
$\sigma _0$ (see $3^0$) corresponds to the first moment 
when a confluence of roots of $P$ or of a root of $P$ 
and a root of $P^{(k)}$ occurs 
(such a confluence occurs at latest for $\sigma =1$ 
because $\dot{y}_1=1$ while the smallest and the greatest roots 
of $P$ remain equal respectively to $0$ and $1$).  

\begin{lm}\label{speed}
1) One can define the speeds $\dot{y}_i$, $i=2,\ldots ,d,$ in a unique way 
so that $\dot{y}_i=\dot{\xi }_i$, 
$i=2,\ldots ,d$. 

2) For these speeds one has $0\leq \dot{y}_i\leq 1$. 

3) The speeds are continuous and bounded on $\bar{\Delta}$ and smooth on 
$\Delta$.
\end{lm}

The lemma is proved after the proof of the theorem. 

\begin{rem} 
The lemma implies property a) 
of the fibration from $1^0$. The absence of stationary points in the 
vectorfield on $\Delta$ results from 
$\dot{y}_i\geq 0$, $\dot{y}_1=1$ which implies that $\dot{a}_1<0$. 
As $\gamma ^{-1}$ is a diffeomorphism, the vectorfield on $T$ has no 
stationary points either.
\end{rem} 

$5^0$. The lemma implies that for $\sigma =\sigma _0$ one or several of the 
following things happen:

-- a root $\xi _{i_0}$ of $P^{(k)}$ which is not a root of $P$ 
becomes equal to a root $y_{j_0}$ 
of $P$ of class A different 
from $y_1$, from the smallest and from the greatest one; for 
$\sigma \in [0,\sigma _0)$ one has $\xi _{i_0}<y_{j_0}$; this is the contrary 
to what happens in {\em i)} from 5) of Section~\ref{dimensionCV};

-- the root $y_1$ becomes equal to a root $\xi _{i_1}$ of $P^{(k)}$ (and 
eventually to $y_{i_1}$ if $y_{i_1}$ is a root of class B); for 
$\sigma \in [0,\sigma _0)$ one has $y_1<\xi _{i_1}$ and $\xi _{i_1}$ is not 
a root of $P$; this is the contrary 
to what happens in {\em i)} or {\em ii)} from 5) of Section~\ref{dimensionCV};

-- the root $y_1$ becomes equal to a root $y_{i_2}$ of class A; for 
$\sigma \in [0,\sigma _0)$ one has $y_1<y_{i_2}$; there might be roots of 
$P^{(k)}$ (and eventually roots of $P$ of class B) 
between $y_1$ and $y_{i_2}$; this is the contrary 
to what happens in {\em iii)} from 5) of Section~\ref{dimensionCV}.  
 
\begin{rems}\label{sigma0}
1) If the CV allows the third possibility (i.e. if the third possibility 
leads to no contradiction 
with condition (\ref{Rolle}) and with Section~\ref{dimensionCV}), then it 
does not allow the second or 
the first one with $j_0=i_2$. Indeed, if the third possibility exists, then 
between $y_1$ and $y_{i_2}$ there must be $\mu -k$ roots of $P^{(k)}$ 
counted with the multiplicities where $\mu$ is the sum of the multiplicities 
of $y_1$, $y_{i_2}$ and of all roots of $P$ (if any) between them; if the 
second possibility exists as well, then for $\sigma =\sigma _0$ there must be 
$\mu -k$ roots of $P^{(k)}$ strictly 
between $y_1$ and $y_{i_2}$ which means that for 
$\sigma <\sigma _0$ there were $\mu -k+1$ of them (one must add the root 
$\xi _{i_1}$) -- a contradiction. In the same way one excludes the first 
possibility with $j_0=i_2$.

2) If the third possibility takes place, then $y_{i_2}$ is 
the first to the right w.r.t. $y_1$ of the roots of class A because these 
roots do not change their positions.

3) Part 1) of these remarks implies that if the CV allows several 
possibilities of the above three types, 
with different possible indices $i_0$, $j_0$, $i_1$, $i_2$ to happen, then 
they can happen  
independently and simultaneously (all of them or any part of them). These 
possibilities can be expressed analytically as conditions (we call them 
{\em equalities} further in the text) of the form 
$y_i=\xi _j$ or $y_{i_1}=y_{i_2}$ for $\sigma =\sigma _0$ while for 
$\sigma <\sigma _0$ there holds $y_i>\xi _j$ or $y_i<\xi _j$ or 
$y_{i_1}>y_{i_2}$.     

4) Property b) of the fibration from $1^0$ follows from 1) -- 3); the CVs of 
the strata from ${\cal U}$ are obtained by replacing certain inequalities 
between roots by the corresponding equalities in the sense of 3) from these 
remarks. 
\end{rems}

$6^0$. Denote by ${\cal U}'$ the set of images under $\gamma$ of 
strata of $\Pi$ (we call these images {\em strata} of ${\cal U}'$) 
whose CVs are obtained from the one 
of $T$ by replacing some or all inequalities by the corresponding equalities, 
see part 4) of Remarks~\ref{sigma0}. 

Consider the vectorfield defined on $\Delta \cup {\cal U}'$ by the conditions 
$\dot{y}_1=-1$ and 
$\dot{y}_i=\dot{\xi }_i$, $i=2,\ldots ,d$. 
On each stratum of ${\cal U}'$, when defining the vectorfield, 
some of the 
multiple roots of $P$ and/or $P^{(k)}$ should be considered as several 
coinciding roots of given multiplicities. What we are doing resembles 
an attempt 
``to revert the phase curves of the already constructed vectorfield on 
$\Delta$'' (and it is the case on $\Delta$) but 
we have not proved yet that each point of each stratum of ${\cal U}'$ is a 
limit point of a phase curve of that vectorfield and that each point of 
${\cal U}'$ belongs to $\bar{\Delta}$. Notice that due to the definition 
of the vectorfield each phase curve stays in $\Delta \cup {\cal U}'$ on 
some time interval.

Each phase curve of the vectorfield defines a family $P_{\sigma}$ of 
polynomials. It is convenient to choose as 
parameter again $\sigma \in [0,\sigma _0]$ where the point of the family 
belongs to ${\cal U}'$ for $\sigma =\sigma _0$.

\begin{lm}\label{A1}
For $\sigma <\sigma _0$ and close to $\sigma _0$ the point of the 
family $P_{\sigma}$ belongs to $\Delta$. 
\end{lm}

The lemma is proved after the proof of Lemma~\ref{speed}. It follows from 
the lemma that ${\cal U}'\subset \bar{\Delta}$. Hence, one can set 
${\cal U}=\gamma ^{-1}({\cal U}')$ and property d) 
of the fibration follows.\\

$7^0$. There remains to be proved that the fibration possesses property c). 
To this end prove first that all strata of dimension $0$ are contractible, 
i.e. connected. Recall that a hyperbolic 
polynomial from a stratum of dimension $0$ has exactly 
two distinct roots of class A -- the smallest and the greatest one 
(see 4) of Section~\ref{dimensionCV}). 

The strata of dimension $0$ whose CVs contain only two multiplicities are 
connected. Indeed, the uniqueness of such monic polynomials up to 
transformations $\gamma$, see $2^0$, is obvious -- 
they equal $x^{m_1}(x-1)^{n-m_1}$. 

Prove the 
uniqueness up to a transformation $\gamma$ of all polynomials defining 
strata of dimension $0$ by induction on $q$ (the number of distinct roots 
of $P$). For $q=2$ the uniqueness is proved above. Denote by $A_i$ 
parts (eventually empty) of the CV which are maximal 
packs of consecutive letters $a$. 

Deduce the uniqueness of the stratum $V$ 
defined by the CV 

\[ \vec{v}=(m_1,A_1,(m_2)_a,A_2,(m_3)_a,A_3,\ldots ,(m_{q-1})_a, 
A_{q-1},m_q)\] 
from the uniqueness of the stratum $W$ defined by the CV 

\[ \vec{w}=(m_1,A_1',(m_2)_a,A_2',(m_3)_a,A_3,\ldots ,
(m_{q-1})_a,A_{q-1},m_q) \] 
We denote again the distinct roots of $P$ by 
$0=y_1<\ldots <y_q=1$ (and we change the indices of the roots $\xi _i$ so 
that on $V$, $\xi _2$, $\ldots$, $\xi _{q-1}$ be equal respectively to $y_2$, 
$\ldots$, $y_{q-1}$). The part $A_1'$ (resp. $A_2'$) contains one letter $a$ 
more than $A_1$ (resp. one letter $a$ less than $A_2$). Eventually $A_1'$ can 
be empty.

To do this construct 
a one-parameter family $P_{\sigma}$ 
(depending on $\sigma \in [0,\sigma _0]$) of 
polynomials joining the two strata 
(for $\sigma =0$ we are on $V$, for $\sigma =\sigma _0$ we are on $W$); 
these polynomials belong to the one-dimensional stratum $Z$ 
defined by the CV 

\[ \vec{z}=(m_1,A_1',m_2,A_2,(m_3)_a,A_3,\ldots ,(m_{q-1})_a, 
A_{q-1},m_q)\] 
For the root $y_2$ one has $\dot{y}_2=1$. One defines 
$\dot{y}_i$, $i=3,\ldots ,q-1$ so that $\dot{\xi}_i=\dot{y}_i$. This 
condition defines them in a unique way (see Lemma~\ref{speed}) 
and there exists a unique 
$\sigma _0>0$ for which one obtains $\vec{w}$ as CV (this follows from the 
uniqueness of $W$ -- the ratio $(y_2-y_1)/(y_2-y_q)=y_2/(y_2-1)$ increases 
strictly with $\sigma$ which implies the uniqueness of $\sigma _0$).

\begin{rem}
One has $P_{\sigma }\in V$ only for 
$\sigma =0$, and for $\sigma >0$ one has $y_2>\xi _2$. This can be proved by 
full analogy with Lemma~\ref{A1}.
\end{rem}

For $\sigma =\sigma _0$ no confluence of roots 
of $P$ or of $P$ and $P^{(k)}$ other than the one 
of $y_2$ with the left most root of $A_2$ can take place. This can be deduced 
by a reasoning similar to the one from part 1) of Remarks~\ref{sigma0}.

On the other hand, one can revert the speeds, i.e. for the polynomial 
defining the CV $\vec{w}$ one can set $\dot{y}_2=-1$, 
$\dot{\xi}_i=\dot{y}_i$, $i=3,\ldots ,q-1$ and deform it continuously into 
a polynomial defining the CV $\vec{v}$; the deformation passes through 
polynomials from the stratum $Z$. This means that the polynomials 
defining the strata $V$ and $W$ can be obtained from the family $P_{\sigma}$. 
The uniqueness of the strata of dimension $0$ is proved.\\ 

$8^0$. Prove the contractibility of the set ${\cal U}$. Each of the strata 
of ${\cal U}$ is defined by a finite number of {\em equalities} (see part 3) 
of Remarks~\ref{sigma0}) which replace inequalities that hold in the CV 
defining the stratum $T$. For each stratum of ${\cal U}$ of dimension 
$p>0$ one can construct a fibration in the same way as this was done for 
$T$ and show that the stratum can be retracted to a finite subset of the 
strata from ${\cal U}$ which are all of dimension $<p$. Hence, ${\cal U}$ can 
be retracted on its only stratum of dimension $0$ (it is defined by all 
equalities). By $7^0$ this stratum is a point. Hence, ${\cal U}$ is 
contractible, 
$T$ as well.~~~~$\Box$\\

{\bf Proof of Lemma~\ref{speed}:}\\ 

$1^0$. Fix the index $i$ of a root of class B. Recall that we denote 
by $m_{\nu}$ the 
multiplicity of the root $y_{\nu}$. Set 
$G_{i,\nu}=(\partial (\xi _i)/\partial (y_{\nu}))$.  
One has 

\[ \dot{\xi}_i=\sum _{\nu =1}^dm_{\nu}G_{i,\nu}\dot{y}_{\nu}~.\]
Hence, the condition $\dot{\xi}_i=\dot{y}_i$ for $i=2,\ldots ,d$ reads:

\begin{equation}\label{system}
\dot{y}_i=\sum _{\nu=1}^dm_{\nu}G_{i,\nu}\dot{y}_{\nu}~~,~~
i=2,\ldots ,d
\end{equation}
Further in the proof ``vector'' means ``$(d-1)$-vector-column''. 
Denote by $V$ the vector with components $\dot{y}_i$. 
Hence, the last system 
can be presented in the form $V=GV+H~(*)$ or $(I-G)V=H$ where $H$ is the 
vector with entries $m_1G_{i,1}$, $2\leq i\leq d$ (recall that 
$\dot{y}_1=1$) and $G$ is the matrix with entries 
$G_{i,\nu}$, $i,\nu =2,\ldots ,d$.\\ 
  
$2^0$. Like in the proof of Proposition~\ref{stratificationprop} 
one shows that the 
matrix $I-G$ is diagonally dominated. Hence, system (\ref{system}) has a 
unique solution $V$. Moreover, its components are all non-negative. Indeed, 
one has $m_1G_{i,1}\geq 0$ for $i=2,\ldots ,d$,  
all entries of the matrix $G$ are non-negative 
(see Lemma~\ref{positive} and Corollary~\ref{positivecor}), and one can 
present $V$ as a convergent series $H+GH+G^2H+\ldots$ whose terms are 
vectors with non-negative entries. This proves 1) and the left 
inequality of 2).\\  

$3^0$. To prove the right inequality of 2) denote by $V_0$ the vector 
whose components are units; write equation $(*)$ in the form 
$(V-V_0)=G(V-V_0)+H+GV_0-V_0$ and observe that all components of the 
vector $H+GV_0-V_0$ are non-positive (this can be deduced from  
Corollary~\ref{estimbis}). Like in $2^0$ we prove that the vector $V-V_0$ 
is with non-positive components. This proves the right inequality of 2).\\ 

$4^0$. Boundedness and continuity of the speeds $\dot{y}_i$ on $\bar{\Delta}$ 
follows from the 
boundedness and continuity of $G$ on $\bar{\Delta}$ (which is compact), and 
from the fact that 
the matrix $I-G$ is uniformly diagonally dominated for any 
point of $\bar{\Delta}$ (see Corollary~\ref{estimbis}). Smoothness of the 
speeds in $\Delta$ follows from the fact that the entries of $G$ are smooth 
there -- all roots $x_j^{(k)}$ are smooth functions of $x_i$ 
inside $\Pi$, i.e. when $x_i$ are distinct.~~~~$\Box$\\ 

{\bf Proof of Lemma~\ref{A1}:}\\ 

$1^0$. We show that for $\sigma <\sigma _0$ and sufficiently close to 
$\sigma _0$ the CV of $P_{\sigma}$ changes -- at least one equality (see 
part 3) of 
Remarks~\ref{sigma0}) is replaced by the corresponding inequality. Hence,  
either the point of the phase curve belongs to $\Delta$ for all 
$\sigma <\sigma _0$ sufficiently close to $\sigma _0$ or it belongs to a 
stratum $S$ of ${\cal U}'$ of higher dimension than the dimension of the 
initial one $S_0$. 
The same reasoning can be applied then to $S$ instead of $S_0$ which will 
lead to the conclusion that the curve cannot stay on $S$ for 
$\sigma \in (\sigma _0-\varepsilon ,\sigma _0]$ for any $\varepsilon >0$ small 
enough. Hence, the curve passes through $\Delta$ for such $\varepsilon$.\\

$2^0$. If for $\sigma =\sigma _0$ there occurs a 
confluence of two roots of $P$  
(w.r.t. $\sigma <\sigma _0$), 
then it is obvious that the CV has changed. So suppose that there occurs a 
confluence of a root $y_{j_0}$ 
of $P$ and of a root $\xi _{i_0}$ of $P^{(k)}$ without a confluence of 
$y_{j_0}$ with another root of $P$. Hence, $y_{j_0}$ is a root of $P$ of 
multiplicity $\leq k-1$.   

By full analogy with Lemma~\ref{speed}, 
one proves that one has $-1\leq \dot{y}_i\leq 0$ for all indices $i$ of roots 
of class B.\\  

$3^0$. Suppose first 
that $j_0=1$.  
Show that one has 
$-1<\dot{\xi}_{i_0}<0$ which implies that the CV has changed 
(because $\dot{y}_1=-1$).  
One has $\dot{\xi}_{i_0}=\sum _{j=1}^{q}
m_j\frac{\partial \xi _{i_0}}{\partial y_j}\dot{y}_j$ with 
$\frac{\partial \xi _{i_0}}{\partial y_j}>0$ for all $j$ 
(see Lemma~\ref{positive1}) 
and $-1\leq \dot{y}_i\leq 0$.
Moreover, one has $\dot{y}_i=0$ for the smallest and for the greatest 
root of $P$. As $\sum _{j=1}^qm_j\frac{\partial \xi _{i_0}}{\partial y_j}=1$ 
(see (\ref{sumi})), one has $-1<\dot{\xi}_{i_0}<0$.\\

$4^0$. If $j_0\neq 1$, then one has $\dot{y}_{j_0}=0$ (because before the 
confluence $y_{j_0}$ has been a root of class A). Like in $3^0$ one shows that 
$-1<\dot{\xi}_{i_0}<0$. Hence, the CV changes again.~~~~$\Box$

Author's address: Universit\'e de Nice -- Sophia Antipolis, 
Laboratoire de Math\'ematiques, Parc Valrose, 06108 Nice Cedex 2, 
France.e-mail: kostov@math.unice.fr

\end{document}